\documentclass{article}

\usepackage{amssymb,amsmath,theorem,euscript}

\newcounter{sec}

\newcounter{punct}

\def\punct{\refstepcounter{punct}{
\arabic{punct}.  }}

\newtheorem{theorem}{Theorem}
\newtheorem{proposition}[theorem]{Proposition}

\newtheorem{lemma}[theorem]{Lemma}

 \def\ov{\overline}

\begin{document}

\def\OO{\mathrm{O}}
\def\Iw{\mathrm{Iw}}
\def\GLO{\mathrm{GLO}}
\def\Coll{\mathrm{Coll}}
\def\kappa{\varkappa}
\def\B{\mathrm B}

\def\R{\mathbb{R}}
\def\C{\mathbb{C}}

\def\la{\langle}
\def\ra{\rangle}

 \def\cA{\mathcal A}
\def\cB{\mathcal B}
\def\cC{\mathcal C}
\def\cD{\mathcal D}
\def\cE{\mathcal E}
\def\cF{\mathcal F}
\def\cG{\mathcal G}
\def\cH{\mathcal H}
\def\cJ{\mathcal J}
\def\cI{\mathcal I}
\def\cK{\mathcal K}
 \def\cL{\mathcal L}
\def\cM{\mathcal M}
\def\cN{\mathcal N}
 \def\cO{\mathcal O}
\def\cP{\mathcal P}
\def\cQ{\mathcal Q}
\def\cR{\mathcal R}
\def\cS{\mathcal S}
\def\cT{\mathcal T}
\def\cU{\mathcal U}
\def\cV{\mathcal V}
 \def\cW{\mathcal W}
\def\cX{\mathcal X}
 \def\cY{\mathcal Y}
 \def\cZ{\mathcal Z}
\def\0{{\ov 0}}
 \def\1{{\ov 1}}
 \def\frA{\mathfrak A}
 \def\frB{\mathfrak B}
\def\frC{\mathfrak C}
\def\frD{\mathfrak D}
\def\frE{\mathfrak E}
\def\frF{\mathfrak F}
\def\frG{\mathfrak G}
\def\frH{\mathfrak H}
\def\frI{\mathfrak I}
 \def\frJ{\mathfrak J}
 \def\frK{\mathfrak K}
 \def\frL{\mathfrak L}
\def\frM{\mathfrak M}
 \def\frN{\mathfrak N} \def\frO{\mathfrak O} \def\frP{\mathfrak P} \def\frQ{\mathfrak Q} \def\frR{\mathfrak R}
 \def\frS{\mathfrak S} \def\frT{\mathfrak T} \def\frU{\mathfrak U} \def\frV{\mathfrak V} \def\frW{\mathfrak W}
 \def\frX{\mathfrak X} \def\frY{\mathfrak Y} \def\frZ{\mathfrak Z} \def\fra{\mathfrak a} \def\frb{\mathfrak b}
 \def\frc{\mathfrak c} \def\frd{\mathfrak d} \def\fre{\mathfrak e} \def\frf{\mathfrak f} \def\frg{\mathfrak g}
 \def\frh{\mathfrak h} \def\fri{\mathfrak i} \def\frj{\mathfrak j} \def\frk{\mathfrak k} \def\frl{\mathfrak l}
 \def\frm{\mathfrak m} \def\frn{\mathfrak n} \def\fro{\mathfrak o} \def\frp{\mathfrak p} \def\frq{\mathfrak q}
 \def\frr{\mathfrak r} \def\frs{\mathfrak s} \def\frt{\mathfrak t} \def\fru{\mathfrak u} \def\frv{\mathfrak v}
 \def\frw{\mathfrak w} \def\frx{\mathfrak x} \def\fry{\mathfrak y} \def\frz{\mathfrak z} \def\frsp{\mathfrak{sp}}
 \def\bfa{\mathbf a} \def\bfb{\mathbf b} \def\bfc{\mathbf c} \def\bfd{\mathbf d} \def\bfe{\mathbf e} \def\bff{\mathbf f}
 \def\bfg{\mathbf g} \def\bfh{\mathbf h} \def\bfi{\mathbf i} \def\bfj{\mathbf j} \def\bfk{\mathbf k} \def\bfl{\mathbf l}
 \def\bfm{\mathbf m} \def\bfn{\mathbf n} \def\bfo{\mathbf o} \def\bfp{\mathbf p} \def\bfq{\mathbf q} \def\bfr{\mathbf r}
 \def\bfs{\mathbf s} \def\bft{\mathbf t} \def\bfu{\mathbf u} \def\bfv{\mathbf v} \def\bfw{\mathbf w} \def\bfx{\mathbf x}
 \def\bfy{\mathbf y} \def\bfz{\mathbf z} \def\bfA{\mathbf A} \def\bfB{\mathbf B} \def\bfC{\mathbf C} \def\bfD{\mathbf D}
 \def\bfE{\mathbf E} \def\bfF{\mathbf F} \def\bfG{\mathbf G} \def\bfH{\mathbf H} \def\bfI{\mathbf I} \def\bfJ{\mathbf J}
 \def\bfK{\mathbf K} \def\bfL{\mathbf L} \def\bfM{\mathbf M} \def\bfN{\mathbf N} \def\bfO{\mathbf O} \def\bfP{\mathbf P}
 \def\bfQ{\mathbf Q} \def\bfR{\mathbf R} \def\bfS{\mathbf S} \def\bfT{\mathbf T} \def\bfU{\mathbf U} \def\bfV{\mathbf V}
 \def\bfW{\mathbf W} \def\bfX{\mathbf X} \def\bfY{\mathbf Y} \def\bfZ{\mathbf Z} \def\bfw{\mathbf w}
 \def\R {{\mathbb R }} \def\C {{\mathbb C }} \def\Z{{\mathbb Z}} \def\H{{\mathbb H}} \def\K{{\mathbb K}}
 \def\N{{\mathbb N}} \def\Q{{\mathbb Q}} \def\A{{\mathbb A}} \def\T{\mathbb T} \def\P{\mathbb P} \def\G{\mathbb G}
 \def\bbA{\mathbb A} \def\bbB{\mathbb B} \def\bbD{\mathbb D} \def\bbE{\mathbb E} \def\bbF{\mathbb F} \def\bbG{\mathbb G}
 \def\bbI{\mathbb I} \def\bbJ{\mathbb J} \def\bbL{\mathbb L} \def\bbM{\mathbb M} \def\bbN{\mathbb N} \def\bbO{\mathbb O}
 \def\bbP{\mathbb P} \def\bbQ{\mathbb Q} \def\bbS{\mathbb S} \def\bbT{\mathbb T} \def\bbU{\mathbb U} \def\bbV{\mathbb V}
 \def\bbW{\mathbb W} \def\bbX{\mathbb X} \def\bbY{\mathbb Y} \def\kappa{\varkappa} \def\epsilon{\varepsilon}
 \def\phi{\varphi} \def\le{\leqslant} \def\ge{\geqslant}

\def\GL{\mathrm {GL}}
\def\PIw{\mathrm {Piw}}
\def\GLB{\mathrm {GLB}}
\def\SL{\mathrm {SL}}
\def\bGL{\mathbf {GL}}
\def\PGL{\mathrm {PGL}}

\def\bGr{\mathbf {Gr}}
\def\Gr{\mathrm {Gr}}
\def\bFl{\mathbf {Fl}}
\def\Fl{\mathrm {Fl}}

\def\St{\mathrm{St}}

 \newcommand{\Dim}{\mathop {\mathrm {Dim}}\nolimits}
  \newcommand{\codim}{\mathop {\mathrm {codim}}\nolimits}
   \newcommand{\im}{\mathop {\mathrm {im}}\nolimits}
\newcommand{\ind}{\mathop {\mathrm {ind}}\nolimits}
\newcommand{\graph}{\mathop {\mathrm {graph}}\nolimits}

\def\F{\mathbb{F}}

\def\sm{\smallskip}

\begin{center}
\Large\bf
On infinite-dimensional limit of the Steinberg representations 

\sc Yury A. Neretin%
\footnote{Supported by the grant FWF, Project P25142}
\end{center}



{\small We present
a construction of the Steinberg representation  admitting an
automatic pass to an infinite-dimensional limit.}

\medskip

Recall, that the representation theory of 
infinite-dimensional classical groups and infinite symmetric groups is a relatively
old  well-developed topic. For infinite-dimensional groups over finite fields a progress appeared comparatively recently,
in \cite{VK}, \cite{VK2}, \cite{GKV} and in \cite{Ner1}, \cite{Ner2}. This note
contains a  construction intermediate between these works.

\smallskip

{\bf\punct Notation.} Denote by $\F_q$ the field with $q$ elements,
let $\F_q^n$ be the coordinate $n$-dimensional linear space with the standard basis $e_j$.
Denote by  $\GL(n)$ the group of
 all invertible  matrices of order $n$. It acts on $\F_q^n$
 by multiplication $x\mapsto xg$ of a row $x\in \F_q^n$ by a matrix $g\in \GL(n)$.

 Denote by $S_n$ the symmetric group. It is generated by the transpositions $\tau_j=(j,j+1)$, the relations are 
 $\tau_j^2=1$, $(\tau_j \tau_{j+1})^3=1$,
 and $\tau_k\tau_j=\tau_j\tau_k$ for $|k-j|\ge 2$. 

By $\ell_2(Z)$ we denote the space of complex functions on
a finite set $Z$ equipped with the $\ell_2$-inner product. 
 
 \sm

{\bf\punct  Schubert cells.} 
  Denote by $\Fl(n)$
 the space of complete flags $\cV$ in $\F_q^n$,
 $$\cV:0=V_0\subset V_1\subset\dots \subset V_n=\F_q^n,\qquad
\dim V_j=j.$$

Fix a standard flag
$$
\cE:\, E_0\subset E_1\subset\dots \subset E_n,
\qquad E_j=\sum_{i\le j} \F_q e_i.
$$
Denote by $\B(n)\subset \GL(n)$ the stabilizer of $\cE$. It consists of lower-triangular
matrices. For each $\sigma\in S_n$ we consider the flag
$$
\cE^\sigma:\,  E_0^\sigma\subset E_1^\sigma\subset\dots \subset E_n^\sigma, \qquad \text{where}
\qquad E^\sigma_j=\F_q e_{\sigma(1)}\oplus\dots\oplus \F_q e_{\sigma(j)}
.
$$
Orbits of $\B(n)$ on $\Fl(n)$ 
are called {\it Schubert cells}, for details, see \cite{Ful}, 10.2.
 They are  enumerated by elements $\sigma\in S_n$: each orbit contains a  unique flag $\cE^\sigma$,
 we denote such orbit by $X^\sigma$.

\sm

{\bf\punct The Steinberg representation.} 
Denote by $\Fl_j(n)$ the space of incomplete  flags 
containing subspaces of all dimensions except $j$. Denote by $\pi_j:\Fl(n)\to\Fl_j(n)$
the map forgetting $V_j$. 
There is a natural map $\Pi_j:\ell_2(\Fl(n))\to \ell_2(\Fl_j(n))$
defined by
$$
\Pi_jf(\cW)=\frac 1{q+1} \sum_{\cV:\pi_j(\cV)=\cW} f(\cW)
$$
In fact, the summation is taken over all flags
$$
W_0\subset\dots \subset W_{j-1}\subset Y\subset W_{j+1}\subset \dots\subset W_n,
$$
such flags are enumerated by subspaces $Y$ satisfying $W_{j-1}\subset Y\subset W_{j+1}$,
or equivalently over set of lines in the 2-dimensional space $W_{j+1}/W_{j-1}$.

\begin{theorem}
 There exists 
 a unique irreducible  representation of $\GL(n)$, which is contained in
$\ell_2(\Fl(n))$ and does not contained in spaces $\ell_2(\Fl_j(n))$.
\end{theorem}

 This representation is called a {\it Steinberg representation},
for its  fascinating properties, see surveys \cite{Hum}, \cite{Ste}.

\sm

{\bf\punct Definition of the Steinberg representation via reproducing kernels.}
We define a function $k(\cV,\cW)$ on $\Fl(n)\times \Fl(n)$ 
as the number of all pairs $(i,j)$, where $i$, $j$ range in $\{0,1,\dots,n-1\}$,
such that
\begin{equation}
\dim V_i\cap W_j=\dim V_{i+1}\cap W_{j+1}
.
\label{eq:dimdim}
\end{equation}

Define the kernel $K(\cdot,\cdot)$ on $\Fl(n)$ by
$$
K(\cV,\cW)=(-q)^{-k(\cV,\cW)}
.
$$
 By definition, the kernel $K(\cdot,\cdot)$ is $\GL(n)$-invariant.
 
 \sm
 
\begin{proposition}
 The function
 $$\kappa(\cW):=k(\cE,\cW)$$
  is constant on
 Schubert cells $X^\sigma$. The value of $\kappa$ on $X^\sigma$ coincides with the number
$I(\sigma)$ of inversions of $\sigma$. The number of points of $X^\sigma$ coincides with
 $q^{k(\cE,\cW)}$.
 \end{proposition}
 
 Notice that $I(\sigma)$
 coincides with the length of a shortest decomposition of $\sigma$ in product
 of the generators  $\tau_j$.
 
 {\sc Proof.} Let us evaluate the number of points of $X^\sigma$.
 Consider an example. Let $n=6$,
 \begin{equation}
 \sigma=
 \begin{pmatrix}
 0&0&0&1&0&0\\
 0&1&0&0&0&0\\
 0&0&0&0&1&0\\
 1&0&0&0&0&0\\
 0&0&0&0&0&1\\
 0&0&1&0&0&0
 \end{pmatrix}
 \label{eq:matritsa1}
 .
 \end{equation}
 Vectors $e_\sigma$ are rows of this matrix. A $\B(n)$-orbit of the collection 
 $\{e_\sigma\}$
 consists  of arbitrary collections of the type
  $$
 \begin{matrix}
(& *&*&*&\circ&0&0&),\\
(& *&\circ&0&0&0&0&),\\
(& *&*&*&*&\circ&0&),\\
(& \circ&0&0&0&0&0&),\\
(& *&*&*&*&*&\circ&),\\
(& *&*&\circ&0&0&0&),
 \end{matrix}
 $$
 where $*$ denotes arbitrary elements of $\F_q$ and $\circ$ nonzero elements.
 We have $\circ$'s on the former positions of units and $*$'s on positions to the left
 of units. Elements of flags (subspaces) are linear combinations $\sum_{j\le k} c_j e_{\sigma(j)}$.
 Replacements
 $$
 e_{\sigma(j)}\to \lambda e_{\sigma(j)}+\sum\nolimits_{i<j} a_i e_{\sigma(i)},
 \quad \lambda\ne 0,
 $$
 do not change the flag.
 Therefore we can get $1$ on the place of $\circ$'s and $0$ under all units.
Thus we get that any flag in $\B(n)$-orbit of $\cE^\sigma$ is generated
by a collection of vectors
  \begin{equation}
 \begin{matrix}
(& *&*&*&1&0&0&),\\
(& *&1&0&0&0&0&),\\
(& *&0&*&0&1&0&),\\
(& 1&0&0&0&0&0&),\\
(& 0&0&*&0&0&1&),\\
(& 0&0&1&0&0&0&).
 \end{matrix}
 \label{eq:3}
 \end{equation}
 Now for each star we have a unit under this star and a unit to the right of the star.
 This pair of units corresponds to an inversion in $\sigma$.
 
\sm 
 
 Next, let us evaluate the  number $(i,j)$ of pairs satisfying
 (\ref{eq:dimdim}). Dimension of $E_i\cap F_j$ is the number of units
 in the left upper $i\times j$ corner of the matrix $\sigma$. 
 The condition  (\ref{eq:dimdim}) means that $i\times j$ and 
 $(i+1)\times (j+1)$-corners contain the same units. Therefore 
 units in $(i+1)$-th row and $(j+1)$-th column are outside the
  $(i+1)\times (j+1)$-corner. Hence we have $*$ on the $(i+1)(j+1)$-th
  place in (\ref{eq:3}).
  \hfill $\square$

\begin{lemma}
 The kernel $K(\cdot,\cdot)$ is positive definite%
 \footnote{I.e., for any collection of points $\cV_i\in\Fl(n)$, we have
 $\det_{i,j}
 \{K(\cV_i,\cV_j)\}\ge 0$.}.
\end{lemma}

Consider the Euclidean space $H_n$ determined by the reproducing kernel%
\footnote{See, e.g., \cite{Ner-gauss}, Section 7.1.}
$K(\cV,\cW)$. 

\begin{lemma}
The representation of $\GL(n)$ in $H_n$ coincides with the Steinberg representation.
\end{lemma}

{\sc Proofs of lemmas.}  By the Frobenius
reciprocity, any subrepresentation in $\ell^2(\Fl(n))$ contains a $\B(n)$-invariant vector.
 Denote by $\eta$ a
 $\B(n)$-invariant function in the Steinberg subrepresentation $\St$ in $\ell_2$. 
 Denote by $\eta[\sigma]$ its value  on a Schubert cell $X^\sigma$.
  By the definition $\eta$ satisfies the equations
 $\Pi_j\eta=0$. 
 It is easy to see that that these equations have the form
 $$
  q\eta[\tau_j\sigma]+\eta[\sigma]=0 \quad\text{if $I(\tau_j \sigma)>I(\sigma)$}
.$$
These recurrence relations have a unique (up to a constant factor) solution,
namely $\eta[\sigma]=(-q)^{-I(\sigma)}$.
This also proves Theorem 1.
 

Denote by $M(\cdot,\cdot)$ the reproducing kernel determining the subspace
$\St$. This means that the functions 
$$\delta_\cV(\cW)=M(\cV,\cW)$$
are contained in $\St$ and for any function $f$ on $\Fl(n)$ we have
$$
f(\cV)=\la f,\delta_\cV\ra_{\ell_2(\Fl(n))}
.
$$
Since $\St$ is $\GL(n)$-invariant, the kernel $M$ is $\GL(n)$-invariant,
$M(g\cV,g\cW)=M(\cV,\cW)$. Since the action of $\GL(n)$ on $\Fl(n)$
is transitive, the kernel 
 is determined by its values for $\cV=\cE$, i.e. by the function
 $\delta_\cE$. Moreover, $\delta_\cE(\cW)$ is $\B(n)$-invariant, and therefore
 $\delta_\cE(\cW)=s\cdot \eta(\cW)$.
 \hfill $\square$

 \sm
 
{\sc Remark.} This construction of the Steinberg representation is a rephrasing of \cite{CIK}, Theorem 10.2.
 
 \sm
 
{\bf\punct Infinite-dimensional limit. Preliminaries.}  Consider the linear space
$L$, whose vectors are two-side sequences%
$$
x=
(\dots,x_{-1}, x_0, x_1 ,\dots).
$$
such that $x_k=0$ for sufficiently large $k$.
We represent operators in $L$ as infinite matrices
$g=g_{ij}$, where $-\infty<i,j<\infty$.
Denote by $\B(2\infty)$ the group  of all infinite matrices
$g$ such that $g_{ij}=0$ for $i<j$ and $g_{ii}$ are invertible
(i.e., we consider all invertible lower-triangular matrices).
Denote by $\GL(2\infty)$ the group of {\it finitary}%
\footnote{This means that $g-1$ has finite number of nonzero matrix elements.}
 invertible matrices. Denote by $\GLB(2\infty)$ the group of matrices
 generated by $\GL(2\infty)$ and $\B(2\infty)$. The group
 $\GLB(2\infty)$ acts in $L$ by transformations $x\mapsto xg$.

Denote by $E_j$, where $-\infty<j<\infty$,  the subspace in $L$  consisting of vectors 
$$(\dots,x_{j-1},x_j,0,0,\dots).$$
Thus we get the {\it standard flag} $\cE$:
$$
\dots\subset E_{-1}\subset E_0\subset E_1\subset\dots
$$


 We define a flag space $\Fl(2\infty)$ as the space of complete flags coinciding with
 the standard flag in all but a finite number of terms. More precisely,
 consider  flags $\cV$ 
 having the following form. Fix $M\le N$. Set
 $V_j=E_j$ if $j\le M$ and $j\ge N$. Consider
 the finite-dimensional space $E_N/E_M$ and a complete flag in $E_N/E_M$,
 $$
 0=F_0\subset F_1\subset\dots \subset F_{N-M}=E_N/E_M
 .
 $$
 For $0\le \alpha\le N-M$, we set $V_{M+\alpha}$
being the preimage of $F_\alpha$ under the projection 
$E_N\to L/E_M$. 
 
 Setting $M=-n$, $N=n$, we get
 that the space $\Fl(2\infty)$ is an inductive limit
 of the chain
  $$
\dots \longrightarrow \Fl(2n+1) \longrightarrow \Fl(2n+3) \longrightarrow\dots
 $$

   The group $\GLB(2\infty)$ acts
 on the space $\Fl(2\infty)$.

 \sm
 
 {\bf\punct The infinite-dimensional limit of Steinberg representations.} 
We define the function $K(\cV,\cW)$  on $\Fl(2\infty)\times \Fl(2\infty)$
as the number of pairs $(i,j)\in \Z\times \Z$ such that
$$
V_{i+1}\cap  W_{j+1}=V_{i}\cap  W_{j}
.
$$
The kernel $K(\cV,\cW)=(-p)^{-k(\cV,\cW)}$ is positive definite on each space $\Fl(2n+1)$ and therefore it is positive-definite
on the inductive limit $\Fl(2\infty)$. We consider the Hilbert  space determined by the reproducing kernel
$K$ and the unitary representation of $\GL(2\infty)$ in this space.

\sm

{\bf\punct Comparison with earlier papers.}
a) Consider the space $L_+$ consisting of sequences
$(x_0,x_1,\dots)$ such that $x_j=0$ for all but a finite number of $j$. 
The same construction
 gives a Steinberg representation obtained in \cite{GKV}.

\smallskip

b) Grassmannians and flags in the space $L$ were considered in 
 \cite{Ner1}. However, the topic of \cite{Ner1}
 is the group of all continuous transformations of (locally compact Abelian group) $L$;
 this group is larger than $\GLB(2\infty)$. Also, \cite{Ner1} treats another space of flags, which have empty intersection with $\Fl(2\infty)$.

{\tt Math.Dept., University of Vienna,

Institute for Theoretical and Experimental Physics, Moscow

Mech. Math. Dept., Moscow State University,

e-mail: neretin(at) mccme.ru

URL:www.mat.univie.ac.at/$\sim$neretin
}

\end{document}